\documentclass[12pt,a4paper]{article}
\usepackage{amsfonts}
\usepackage{amsmath}
\usepackage{mathrsfs}
\usepackage{graphicx}
\usepackage{amsmath,amsthm,amssymb,amscd}

\renewcommand{\baselinestretch}{1.05}

\newtheorem{thm}{Theorem}[section]

\newtheorem{defn}[thm] {Definition}
\newtheorem{ex}  [thm]{Example}

\theoremstyle{remark}
\newtheorem{rem} [thm]{Remark}

\theoremstyle{claim}

\begin{document}
\title{\bf Systems with Almost Specification Property May Have Zero Entropy}
\author{{Yiwei Dong}\\
{\em\small  School of Mathematical Science,  Fudan University}\\
{\em\small Shanghai 200433, People's Republic of China}\\
{\small   email: dongyiwei0@sina.com }}
\date{}
\maketitle

\renewcommand{\baselinestretch}{1.2}
\large\normalsize

\footnotetext {AMS Review:   37B05; 37B40.  }

\begin{abstract}
It is shown that there exist systems having almost specification property and zero entropy. Since Sigmund has shown that systems with specification property must have positive entropy, this result reveals further the difference between almost specification and specification. Moreover, one can step on to obtain a both sufficient and necessary condition to ensure positive entropy.
\end{abstract}

\section{Introduction}
By a dynamical system $(X,d,T)$, we mean a compact metric space $X$ with metric $d$ and a continuous map $T:X\to X$. We denote by $\mathcal{M}_T(X)$ the $T$-invariant probability measures on $X$. We also denote by $h_{top}(T)$ the topological entropy of $(X,d,T)$ and by $h_{\mu}(T)$ the metric entropy of $\mu\in \mathcal{M}_T(X)$. In this paper, we focus on systems with almost specification property which derives from the classical notion of specification property. So let us recall some developments.

The specification property was first introduced by Bowen in \cite{Bowen1971}. A continuous map $T:X\to X$ satisfies the specification property if for all $\epsilon>0$,  there exists an integer $m(\epsilon)$ such that for any collection $\{I_{j}=[a_{j},b_{j}]\subset\mathbb{N}:j=1,\cdots,k\}$ of finite intervals with $a_{j+1}-b_{j}\geq m(\epsilon)$ for $j=1,\cdots,k-1$ and any $x_{1},\cdots,x_{k}$ in $M$, there exists a point $x\in $ such that
$$d(T^{a_{j}+t}x,T^{t}x_{j})<\epsilon $$ for all $t=0,\cdots,b_{j}-a_{j}$ and $j=1,\cdots,k$. Furthermore, we say $T$ satisfies Bowen's periodic specification if the point $x$ can be chosen to be periodic with period $p$ for every $p\geq b_{k}-a_{1}+m(\epsilon)$.

This seemly strong condition fits well in several systems. For example, the full shift on symbolic systems $X=\Pi_{-\infty}^{+\infty}S$ satisfies specification property, where $S$ is a compact metric space. In particular, if $S$ is discrete space, we obtain the shift on a finite alphabet. In this case, the subshifts of finite type also have the specification property \cite{Sigmund1974}. Besides symbolic systems, one is encouraged to find examples from Axiom A diffeomorphisms. For instance, if $X$ is a $n$-dimensional torus and $T$ is an automorphism induced by a matrix from $SL(n,\mathbb{Z})$ whose eigenvalues are off the unit circle, then $(X,T)$ has the specification property. Meanwhile, it was shown by Blokh that any continuous topologically mixing interval map has the
specification property \cite{Blokh,Buzzi}.
 Beyond these discrete cases, it is easy to define the analogue of specification property for one parameter flow $\{T_t:X\to X\}_{t\in \mathbb{R}}$. One good example is the C-dense Axiom A flows \cite{Bowen1972}. In particular, geodesic flows on manifolds
of negative curvature belong to this class \cite{Bowen1972}.

Since its introduction, systems with specification property have been extensively studied. Here we list some known results for such systems. From the viewpoint of ergodic theory, Sigmund \cite{Sigmund1974} gave a description. He showed that for any non-empty, compact and
connected subset $V\subset \mathcal {M}_{T}(X)$, there exists a dense subset $Y\subset X$ such that the time average of atom measure along the orbit of $x\in Y$ diverges everywhere to $V$. Recently, Climenhaga and Thompson \cite{Climenhaga-Thompson2013} studied equilibrium states beyond specification and the Bowen property. They also \cite{Climenhaga-Thompson2014} introduced the notions of obstructions
to expansivity and specification, and showed that if the entropy of such obstructions is smaller
than the topological entropy of the map, then there is a unique measure of maximal entropy.

Moreover, from the viewpoint of dynamical complexity, specification property has been proved to be a indicator of fairly strong chaos in one way or another. For example, Sigmund \cite{DGS} showed that systems with specification property have positive entropy, which indicates that they have topological chaos. Furthermore, since positive entropy implies Li-Yorke chaos \cite{Forti}, one gets that systems with specification property have Li-Yorke chaos. On the other hand, Sklar and Smital \cite{Sklar-Smital} showed that systems with specification property exhibit distributional chaos. Moreover, if there is a pair of distal points, the distributional chaos is satisfied in a strong sense \cite{Oprocha-Stefankova}. Meanwhile, Oprocha posed a question in \cite{Oprocha-Stefankova} concerning the existence of invariant distributionally
chaotic scrambled sets. Recently, two progresses have been obtained toward this direction under specification property and some other conditions which we refer to \cite{Dolezelova,Wang-Wang} for details. In addition, Oprocha \emph{et al} \cite{GKLOP} defined and used two new versions of specification properties to study the chaos on hyperspaces.

Besides, if one considers from the dimensional viewpoint, there are also rich results. For example, Takens and Verbitskiy studied the multifractal spectrum of local entropies in \cite{Takens-Verbitskiy1999}. They also \cite{Takens-Verbitskiy} studied the entropy spectrum for Birkhoff averages which was later extended to Banach valued Birkhoff ergodic averages \cite{FLP}. Meanwhile, Thompson studied the pressure of level sets (which generalizes entropy) in \cite{Thompson2009}. In particular, it is shown that the irregular points have full entropy \cite{CTS} and full pressure \cite{Thompson2010}. Similar results also hold for Hausdorff dimension and shifts with specification property \cite{Barreira-Saussol2000}.

Finally, from the topological viewpoint, Li and Wu \cite{Li-Wu2014} proved that the set of irregular points is either is empty or residual.

During the evolvement of specification property, a natural but important generalization appeared in the study of large deviation \cite{PS2005}. It was called $g$-almost product property and renamed now as almost specification property \cite{Thompson2012}. The only difference between them is that in the latter case, the mistake function $g$ can depend on $\epsilon$, and thus is a priori slightly weaker than the former.

Recently, systems with almost specification property have gained more and attentions. For example, Pfsiter and Sullivan provided a variational principle for saturated sets in \cite{PS2007}. Moreover, under the additional condition of expansiveness, Yamamoto \cite{Yamamoto} gave the topological pressure formula for periodic orbits. Later on, for continuous function $\varphi:X\to \mathbb{R}$, Thompson \cite{Thompson2012} studied the set $I_{\varphi}$ of $\varphi$-irregular points and showed that $I_{\varphi}$ either is empty or has full entropy. Thompson \cite{Climenhaga-Thompson2012} also adapt almost specification  to the study of symbolic spaces with a non-uniform
structure and get some results concerning intrinsic ergodicity. Meanwhile, under the additional condition of uniform separation which was introduced in \cite{PS2007}, Zhou and Chen \cite{Zhou-Chen} divided the historic set into different level sets and used
topological pressure to describe the size of these level sets. Besides, Oprocha \emph{et al} \cite{KKO} compared almost specification and other similar notions such as asymptotic average shadowing and average shadowing.

Carefully readers may find that many properties can be inherited from specification to almost specification. It seems that there is little difference between these two definitions. However, in the system of $\beta$-shift, they behave quite differently. It is well known that every $\beta$-shift has almost specification property \cite{PS2005}, while the set of $\beta$ for which the $\beta$-shift has specification property has zero Lebesgue measure \cite{Buzzi,Schmeling}. Thus almost specification is a nontrivial generalization of specification.

Bearing this in mind, one can still ask for other characters to distinguish these two notions. This is exactly what we are focusing in this paper. More precisely, we show that there exist systems $(X,d,T)$ satisfying almost specification property and that $h_{top}(T)=0$, in contrasting with Sigmund's result that systems with specification property must have positive entropy \cite[Proposition 21.6]{DGS}. To be concrete, we give a sufficient condition to ensure that $(X,d,T)$ has the required properties.

%Since almost specification property is slightly more general than $g$-almost product property \cite{Thompson2012}, we restrict our attention to systems with almost specification property in this paper.
\begin{thm}\label{counter-example}
Let $(X,d,T)$ be a dynamical system. If there exists a fixed point and $N\in \mathbb{N}$ such that for any $x,y\in X$, there exists $0\leq i\leq N-1$ so that $T^i(x)=T^i(y)$, then $(X,d,T)$ has almost specification property and zero entropy.
\end{thm}

Based on this observation, one is motivated to find the sufficient condition to ensure that $(X,d,T)$ has positive entropy. In fact, one can step further to have a both sufficient and necessary condition to ensure $h_{top}(T)>0$. More precisely, we have the following theorem.

\begin{thm}\label{main-thm}
Let $(X,d,T)$ be a dynamical system satisfying almost specification property with mistake function $g$. Then $h_{top}(T)>0$ if and only if there exist $\sigma>2\delta>0$ and $x,y\in X$ such that $x,y$ are $(2g;N,\sigma)$-separated for some $N\in \mathbb{N}$ with $k_g(\delta)\leq N$.
\end{thm}

\section{Preliminaries}

\subsection{The almost specification property}

\begin{defn}
Let $\epsilon_0>0$. A function $g:\mathbb{N}\times (0,\epsilon_0)\to \mathbb{N}$ is called a mistake function if for all $\epsilon\in (0,\varepsilon_0)$ and all $n\in \mathbb{N}$, $g(n,\epsilon)\leq g(n+1,\epsilon)$ and
$$\lim_{n}\frac{g(n,\epsilon)}{n}=0.$$
Given a mistake function $g$, if $\epsilon\geq \epsilon_0$, we define $g(n,\epsilon)=g(n,\epsilon_0)$.
\end{defn}
\begin{defn}
Let $g$ be a mistake function and $\epsilon>0$. For $n\in \mathbb{N}$ large enough such that $g(n,\epsilon)<\epsilon$, we define
$$I(g;n,\epsilon):=\{\Lambda\subset\{0,1,\cdots,n-1\}:\#\Lambda\geq n-g(n,\epsilon)\},$$
where $\#\Lambda$ denotes the cardinality of $\Lambda$.
\end{defn}

\begin{defn}\label{Bowen-ball}
For a finite set of indices $\Lambda\subset\{0,1,\cdots,n-1\}$, we define the Bowen distance of $x,y\in X$ along $\Lambda$ by
$$d_{\Lambda}(x,y):=\max_{j\in \Lambda}\{d(T^jx,T^jy)\}$$
and the Bowen ball of radius $\epsilon$ centered at $x$ by
$$B_{\Lambda}(x,\epsilon):=\{y\in X:d_{\Lambda}(x,y)<\epsilon\}.$$
When $g(n,\epsilon)<n$, we define the $(g;n,\epsilon)$ Bowen ball centered at $x$ as
$$B_n(g;x,\epsilon):=\{y\in X:y\in B_{\Lambda}(x,\epsilon)~\textrm{for some}~\Lambda\in I(g;n,\epsilon)\}=\bigcup_{\Lambda\in I(g;n,\epsilon)}B_{\Lambda}(x,\epsilon).$$

\end{defn}
Now we are in a position to define almost specification property.

\begin{defn}\label{def-of-alm-spe}
The dynamical system $(X,d,T)$ has almost specification property with mistake function $g$, if for any
$\epsilon_{1},\cdots,\epsilon_{m}>0$, there exist integers $k_g(\epsilon_1),\cdots,k_g(\epsilon_m)$ such that for any points $x_{1},\cdots,x_{m}\in X$, and integers $n_{1}\geq k_g(\epsilon_{1}),\cdots,n_{m}\geq k_g(\epsilon_{m})$, we can find a point $z\in X$ such that
\begin{equation*}
  T^{l_{j}}(z)\in B_{n_{j}}(g;x_{j},\epsilon_{j}),~j=1,\cdots,m,
\end{equation*}
where $n_{0}=0~\textrm{and}~l_{j}=\sum_{s=0}^{j-1}n_{s}$.

\end{defn}
In other words, the appropriate part of the orbit of $z$ $\epsilon_{j}$-traces
the orbit of $x_{j}$ with at most $g(n_{j},\epsilon_{j})$ mistakes. However, in the case of specification, the appropriate part of the orbit of $z$ $\epsilon$-traces
the orbit of $x_{j}$ with at most $m(\epsilon)$ mistakes. Thus if one lets $g(n,\epsilon)=m(\epsilon)$ and $k_g(\epsilon)=m(\epsilon)+1$, it is easy to see that definition \ref{def-of-alm-spe} is satisfied for systems with specification property and $\epsilon_{1}=\cdots=\epsilon_{m}=\epsilon$. Finally, using a trick to replace $\epsilon$ by $\epsilon_{1},\cdots,\epsilon_{m}$, Pfister and Sullivan \cite{PS2007} showed that the specification property implies the almost
specification property.
\subsection{Topological entropy}

For $\epsilon>0$, two points $x$ and $y$ are $(n,\epsilon)$-separated if
$$\max_{0\leq j\leq n-1}\{d(T^{j}(x),T^{j}(y))\}>\epsilon.$$
For $\delta>0$ and $\epsilon>0$, two points $x$ and $y$ are $(g;n,\epsilon)$-separated if
\begin{equation}\label{separated}
  \#\{j:d(T^{j}(x),T^{j}(y))>\epsilon,~0\leq j\leq n-1\}> g(n,\epsilon).
\end{equation}
A subset $Z\subset X$ is $(g;n,\epsilon)$-separated if every two points $x,y\in Z$ are $(g;n,\epsilon)$-separated. Let $Y\subset X$ and we define
\begin{eqnarray*}
  s_n(Y,\epsilon)&:=&~\sup\{\#Z:Z\subset Y~\textrm{is} ~(n,\epsilon)-\textrm{separated}\}, \\
  s_n(g;Y,\epsilon)&:=&~\sup\{\#Z:Z\subset Y~\textrm{is} ~(g;n,\epsilon)-\textrm{separated}\}.
\end{eqnarray*}
We have the following Bowen's definition of topological entropy
\begin{equation}\label{Bowen-entropy}
  h_{top}(T)=\lim_{\epsilon\to 0}\lim_{n\to \infty}\frac{\log s_n(X,\epsilon)}{n}.
\end{equation}
We also have the following equivalent definition of topological entropy \cite[Theorem 3.6]{Thompson2012}
\begin{equation}\label{Thompson-entropy}
  h_{top}(T)=\lim_{\epsilon\to 0}\lim_{n\to \infty}\frac{\log s_n(g;X,\epsilon)}{n}.
\end{equation}

\section{Proof of theorem \ref{counter-example} and some examples}
\textbf{Proof of theorem \ref{counter-example}}. Suppose $z\in X$ is the fixed point. Fix any $x,y\in X$, then there exists an $0\leq i\leq N-1$ with
$T^i(x)=T^i(y)$.
One sees immediately that $T^j(x)=T^j(y)$ for any $j\geq i$. In particular, for any $x,y\in X$, we have $T^{N-1}(x)=T^{N-1}(y)$
and thus
\begin{equation}\label{N}
  T^{j}(x)=T^{j}(y),~j\geq N-1.
\end{equation}
\begin{itemize}
  \item By (\ref{N}), the fixed point $z=f^{N-1}(z)=f^{N-1}(x)$ for any $x\in X$. Thus
  \begin{equation}\label{N-x}
    T^{j}(x)=z,~j\geq N-1 ~~\textrm{for any}~~x\in X.
  \end{equation}
   Now let $g\equiv k_g\equiv N$, then for any
$\epsilon_{1},\cdots,\epsilon_{m}>0$, any points $x_{1},\cdots,x_{m}\in X$, and integers $n_{1}\geq N,\cdots,n_{m}\geq N$, one sees that
\begin{equation}\label{almost}
  T^{l_{j}}(z)\in B_{n_{j}}(g;x_{j},\epsilon_{j}),~j=1,\cdots,m,
\end{equation}
where $n_{0}=0~\textrm{and}~l_{j}=\sum_{s=0}^{j-1}n_{s}$. Thus $(X,d,T)$ satisfies almost specification property.
  \item On the other hand, one sees from (\ref{N-x}) that for any $x\in X$ and $\epsilon>0$,
$$B_n(x,\epsilon)=B_N(x,\epsilon),~n\geq N$$
which implies that
$$s_n(X,\epsilon)=s_N(X,\epsilon),~n\geq N.$$
Therefore,
$$h_{top}(T)=\lim_{\epsilon\to 0}\lim_{n\to \infty}\frac{\log s_n(X,\epsilon)}{n}=\lim_{\epsilon\to 0}\lim_{n\to \infty}\frac{\log s_N(X,\epsilon)}{n}=0.$$
\end{itemize}
\qed

The following simple example gives some intuition.
\begin{ex}Let $X=\{a,b\}$ and $T:X\to X$ be
$$T(a)=b, ~T(b)=b.$$
\end{ex}
One easily checks that $T:X\to X$ is continuous and satisfies almost specification property with mistake function $g\equiv 1$. However, since $X$ is finite, $h_{top}(X)=0$.

Let us give another non-trivial example in the spirit of the condition in theorem \ref{counter-example}.
\begin{ex}\label{non-trivial}
Let $S$ be a compact metric space with metric $\widetilde{d}$. Let $\widetilde{X}=\prod_{0}^{+\infty}S$ whose element is a unilateral sequence $\underline{x}=(x_0,x_1,x_2,\cdots)$ with $x_i\in S,~ i\geq 0$. The distance between two points $\underline{x},\underline{y}\in \widetilde{X}$ is defined as
$$d(\underline{x},\underline{y}):=\sum_{i=0}^{+\infty}\frac{\widetilde{d}(x_i,y_i)}{2^i}.$$
The left shift $T:\widetilde{X}\to \widetilde{X}$ is defined as $(T\underline{x})_n=x_{n+1}$, $n\geq 0$. Select an arbitrary $w\in S$. Let us consider $X\subset \widetilde{X}$ whose element is $\underline{y}=(y_0,y_1,\cdots,y_{N-1},w,w,w,\cdots)$, i.e. $y_j=w$, $j\geq N$. One easily sees that $X$ is a closed and $T$-invariant subset of $\widetilde{X}$. Thus one can consider subsystem $T:X\to X$. It is not hard to see that $(X,d,T)$ satisfies the condition of theorem \ref{counter-example}, i.e. for any $\underline{x},\underline{y}\in X$, there exist an $0\leq i\leq N-1$ such that
$T^i(\underline{x})=T^i(\underline{y}).$
\end{ex}
\begin{rem}
We mention here that the condition of the existence of a fixed point should not be omited. To see this, we assume that there exists no fixed point. Then by (\ref{N}), there exists a common point $w\in X$ such that $w=T^{N-1}(x)=T^{N-1}(w)$ for any $x\in X$. Let
$$\epsilon=\min\{d(T^i(w),T^j(w)):0\leq i<j\leq N-1\}>0.$$
Then for any mistake function $g$ and any function $k_g:(0,+\infty)\to \mathbb{N}$, if one selects
$$n_1=l(N-1)\geq k_g(\epsilon)~\textrm{for some}~l\in\mathbb{N}~~\textrm{and}~~n_2= k_g(\epsilon),$$
it is not hard to see that there exists no $z\in X$ such that
$$z\in B_{n_1}(g;w,\epsilon)~~\textrm{and}~~T^{n_1}(z)\in B_{n_2}(g;Tw,\epsilon)$$
since $T^{N-1}(z)\equiv w$.

\end{rem}

\section{Proof of Theorem \ref{main-thm}}
Since definition (\ref{Thompson-entropy}) is crucial and there is no explicit proof in \cite{Thompson2012}, here we give a version for completeness.

\textbf{Proof of (\ref{Thompson-entropy})}. Consider two cases.
\begin{itemize}
  \item $h_{top}(T)=0$. By definition, $s_n(g;X,\epsilon)\leq s_n(X,\epsilon)$. Then one sees that
  $$0=h_{top}(T)=\lim_{\epsilon\to 0}\lim_{n\to \infty}\frac{\log s_n(X,\epsilon)}{n}\geq\lim_{\epsilon\to 0}\lim_{n\to \infty}\frac{\log s_n(g;X,\epsilon)}{n}\geq0.$$
  Thus (\ref{Thompson-entropy}) holds.
  \item $h_{top}(T)=h>0$. Since $s_n(g;X,\epsilon)$ is non-increasing in $\epsilon$, it is enough to prove that for any $\eta>0$, there exist $\epsilon>0$ and $M\in \mathbb{N}$ such that for any $n\geq M$,
      $$s_n(g;X,\epsilon)\geq exp(n(h-\eta)).$$
      Now we fix any $\eta>0$, by variational principle \cite{Walters}, there exists an ergodic measure $\mu \in \mathcal{M}_T(X)$ such that
      $$h_{\mu}(T)\geq h-\frac{\eta}{2}.$$
      Let us recall the following modified Katok entropy formula \cite[Theorem 3.4]{Thompson2012} for ergodic $\mu\in \mathcal{M}_T(X)$ and $\gamma\in(0,1)$:
      $$h_{\mu}(T)=\lim_{\epsilon\to 0}\lim_{n\to \infty}\frac{1}{n}\log(inf \{ s_n(g;Y,\epsilon):Y\subset X,~\mu(Y)\geq 1-\gamma\}).$$
      In particular, we have $\mu(X)=1>1-\gamma$. Thus there exist $\epsilon>0$ and $M\in \mathbb{N}$ such that for any $n\geq M$,
      $$s_n(g;X,\epsilon)\geq exp(n(h_{\mu}(T)-\frac{\eta}{2}))\geq exp(n(h-\eta)).$$
\end{itemize}\qed

Now we set our proof of theorem \ref{main-thm} in two parts.

\textbf{The if part}. For any $m\in \mathbb{N}$, let $\Sigma_2^m=\{x,y\}^m$ whose element is $\xi_m=(w_1,w_2,\cdots, w_m)$ such that $w_i\in \{x,y\}$, $1\leq i\leq m$.
For $j\in \mathbb{N}$, let $l_j=(j-1)N$. Since $N\geq k_g(\delta)$, by the almost specification property, for any $\xi_m\in \Sigma_2^m$, we can find a point $z_{\xi_m}\in X$ such that
$$T^{l_j}(z_{\xi_m})\in B_N(g;w_j,\delta),~j=1,\cdots,m.$$
Because $x,y ~\textrm{are}~(2g;N,\sigma)-\textrm{separated},$  by (\ref{separated}) and definition \ref{Bowen-ball}, \ref{def-of-alm-spe}, one sees that if $\xi_m\neq \xi_m'$, then $z_{\xi_m}$ and $z_{\xi_m'}$ are $(mN,\sigma-2\delta)$-separated. Thus one has
$$s_{mN}(X,\sigma-2\delta)\geq 2^m,$$
which implies that
$$h_{top}(f)=\lim_{\epsilon\to 0}\lim_{n\to \infty}\frac{\log s_n(X,\epsilon)}{n}\geq\lim_{m\to \infty}\frac{\log s_{mN}(X,\sigma-2\delta)}{mN}\geq\lim_{m\to \infty}\frac{\log(2^m)}{mN}=\frac{\log2}{N}.$$
Here, the second inequality comes from the fact that $s_n(X,\epsilon)$ is non-increasing in $\epsilon$.

\textbf{The only if part}. Suppose $h_{top}(T)=h>0$. Let $\widetilde{g}=2g$, then it is easy to see that
$$B_{n}(g;x,\epsilon)\subset B_{n}(\widetilde{g};x,\epsilon)$$
for any $n\in \mathbb{N}$, $x\in X$ and $\epsilon>0$. Therefore, after checking definition \ref{def-of-alm-spe}, one sees that $(X,d,T)$ also satisfies almost specification property with mistake function $\widetilde{g}$ and $k_{\widetilde{g}}=k_g$.
Thus by the definition in (\ref{Thompson-entropy}) with $\widetilde{g}$ instead of $g$, there exist $\sigma>0$ and $N_1\geq\frac{2\log2}{h}$ such that for any $n\geq N_1$,
\begin{equation}\label{separation}
  s_n(\widetilde{g};X,\sigma)\geq \exp(n\cdot\frac{h}{2})\geq2.
\end{equation}
If we choose $\delta=\frac{\sigma}{3}$ and $N=\max\{N_1,k_g(\frac{\sigma}{3})\}$, then by (\ref{separation}), there exist $x,y\in X$ such that $x,y$ are $(2g;N,\sigma)$-separated.\qed

\end{document}